\documentclass{amsart}
 
\usepackage{amssymb,amsmath,amsfonts,amscd} 
\usepackage{psfrag,graphicx} 
\usepackage[all,dvips,cmtip]{xy}

\title[Crepant resolutions of weighted projective 
spaces]{Crepant resolutions of weighted projective spaces and quantum deformations} 
\author{Samuel Boissi\`ere \and \'Etienne Mann \and  
Fabio Perroni}\thanks{F.P. was partially supported by SNF, No 200020-107464/1} 
\date{\today} 
\address{Samuel Boissi{\`e}re, Laboratoire J.A.Dieudonn\'e UMR CNRS 6621, 
         Universit\'e de Nice Sophia-Antipolis, Parc Valrose, 06108 Nice} 
\email{sb@math.unice.fr} 
\address{\'Etienne Mann, SISSA, Via Beirut 2-4, 34014 Trieste, Italy} 
\email{mann@sissa.it} 
\address{Fabio Perroni, Institut f\"ur Mathematik, Universit\"at Z\"urich, Winterthurerstrasse 190, 
8057 Z\"urich, Switzerland} 
\email{fabio.perroni@math.unizh.ch}

\def\cf{\textit{cf.}\kern.3em} 
 
\def\ie{\textit{i.e.} } 
 
\def\resp{\textit{resp.}\kern.3em} 
\renewcommand{\k}{\kern2pt} 
 
 \numberwithin{equation}{section} \makeatletter 
 
\newcommand{\IR}{\mathbb{R}} 
\newcommand{\IN}{\mathbb{N}} 
\newcommand{\IZ}{\mathbb{Z}} 
\newcommand{\IQ}{\mathbb{Q}} 
\newcommand{\IC}{\mathbb{C}} 
 
\newcommand{\IP}{\mathbb{P}} 
\newcommand{\IV}{\mathbb{V}} 
 
\newcommand{\ii}{\mathrm{i}}

\newcommand{\cF}{\mathcal F}
\newcommand{\cM}{\mathcal M} 
\newcommand{\cO}{\mathcal O}
\newcommand{\cR}{\mathcal R}
 
\newcommand{\cV}{\mathcal V} 
\newcommand{\cX}{\mathcal X}

 \DeclareMathOperator{\CR}{\mathrm CR}

 \DeclareMathOperator{\PD}{\mathrm PD}  
 
\theoremstyle{plain} 
\newtheorem{theorem}[equation]{Theorem} 
\newtheorem{proposition}[equation]{Proposition} 
\newtheorem{lemma}[equation]{Lemma}

\theoremstyle{definition} 
\newtheorem{example}[equation]{Example} 
\newtheorem{remark}[equation]{Remark} 
\newtheorem{notation}[equation]{Notation}

\newcommand{\al}{\alpha} 
\newcommand{\be}{\beta} 
\newcommand{\cd}{\cdot} 
 
\newcommand{\De}{\Delta} 
\newcommand{\de}{\delta} 
\newcommand{\e}{\epsilon} 
\newcommand{\fr}{\frac } 
\newcommand{\ga}{\gamma} 
\newcommand{\Ga}{\Gamma}

\newcommand{\Lam}{\Lambda} 
\newcommand{\lan}{\langle}

\newcommand{\mcal}{\mathcal} 
 
\newcommand{\mb}{\mbox} 
 
\newcommand{\nf}{\normalfont}

\newcommand{\op}{\oplus} 
\newcommand{\ot}{\otimes} 
 
\newcommand{\ran}{\rangle} 
\newcommand{\ra}{\rightarrow} 
\newcommand{\si}{\sigma} 
\newcommand{\Si}{\Sigma} 
\newcommand{\ti}{\tilde} 
\newcommand{\ul}{\underline}

\numberwithin{equation}{section} 
  
\DeclareMathOperator{\Spec}{Spec} 
 
\DeclareMathOperator{\Graph}{Graph} 
 
\begin{document}

\begin{abstract} 
We compare the Chen-Ruan cohomology ring of the weighted
projective spaces $\IP(1,3,4,4)$ and $\IP(1,...,1,n)$
with the cohomology ring of their crepant resolutions.\\
In both cases, we prove that the Chen-Ruan cohomology
ring is isomorphic to the quantum corrected cohomology 
ring of the crepant resolution after suitable evaluation 
of the quantum parameters.
For this, we prove a formula for the Gromov-Witten invariants
of the resolution of a transversal ${\rm A}_3$ singularity.
\end{abstract} 
\maketitle 
\section{Introduction} 

Given a complex orbifold $\cX$, Chen and Ruan defined the so called 
Chen-Ruan cohomology ring of $\cX$ \cite{CRnco}, it is denoted by $H^\star_{\rm CR}(\cX)$
(see \cite{AGVaoqp} and \cite{AGVgwdms} for the definition in the algebraic case).
The \textit{Cohomological Crepant Resolution Conjecture}, as proposed by Y. 
Ruan \cite{Ru}, predicts the existence of an isomorphism between the 
Chen-Ruan cohomology ring of a Gorenstein orbifold $ {\cX}$ and the 
quantum corrected cohomology ring of any crepant resolution $\rho 
:Z\rightarrow  X$ of the coarse moduli space $X$ of $ {\cX}$.  
The quantum corrected cohomology ring of $Z$ is the ring obtained 
from the small quantum cohomology of $Z$ after specialization
of some quantum parameters and is denoted by $H^\star_\rho(Z)(q_1,...,q_m)$ 
(see Sec. \ref{sec:CCRC}). The 
conjecture belongs to the so called generalized McKay correspondence. 

The main examples used to study the conjecture are the following: 
\begin{itemize} 
\item the Hilbert scheme of $r$ points on a projective surface $S$ 
  (see \cite{LQ} for $r=2$; \cite{ELQ} and \cite{LL} for $S=\IP^2$ $r=3$; 
  \cite{FG},  \cite{U}  and \cite{QW} for $r$ 
  general and $S$ with numerically trivial canonical class;
we remark that the proofs of \cite{FG} and \cite{U} use
the computation of $H^\star({\rm Hilb}^r(S))$ done in \cite{LS2},
while \cite{QW} is self-contained); 
\item the Hilbert scheme of $r$ points on a quasi-projective surface 
  $S$ carrying a holomorphic symplectic form (see \cite{LQW}, 
  \cite{LS} and \cite{V}); 
\item  the quotient $V/G$, where $V$ is a complex symplectic vector space and $G$ is a 
  finite subgroup of $\mathrm{Sp}(V)$ (see \cite{GK});
\item orbifolds with transversal ${\rm A}_n$ singularities (see \cite{P});
\item wreath product orbifolds (see \cite{Tomoo}).
\end{itemize} 
 
The main result of the present paper 
concerns $\IP(1,3,4,4)$.  The singular locus of its coarse 
moduli space $|\IP(1,3,4,4)|$ is the disjoint 
union of an isolated singularity of type $\fr{1}{3}(1,1,1)$ (we use 
Reid's notation \cite{YPG}) and a transversal ${\rm A}_3$ singularity. 
The quantum corrected cohomology ring of the crepant resolution 
$\rho:Z\ra X$
has four quantum parameters: $q_1,q_2,q_3$ and $q_4$; 
the parameters $q_1,q_2,q_3$ 
come from the resolution of the transversal singularity and 
$q_4$ comes from the isolated 
singularity.  We compute the $3$-points genus zero 
Gromov-Witten invariants of curves which are contracted by $\rho$
and which are contained in the part of the exceptional divisor 
that contracts to the transversal singularity (Theorem \ref{GW1344}).
Then we prove that, for $(q_1,q_2,q_3,q_4) \in \{ (\ii,\ii,\ii,0), (-\ii,-\ii,-\ii,0) \}$ 
  there is a ring isomorphism (Theorem \ref{th2})
\begin{equation*} 
H^{\star}_{\rho}(Z;\IC)(q_1,q_2,q_3,q_4) \cong H^{\star}_{\rm CR}(\IP(1,3,4,4);\IC). 
\end{equation*}

\medskip
In the case of $\IP(1,...,1,n)$, we prove that the Chen-Ruan
cohomology ring is isomorphic to the cohomology ring of
the crepant resolution $\rho :Z \ra |\IP(1,...,1,n)|$ (Proposition \ref{p11n}). 
In this case we will always assume that the dimension
of $\IP(1,...,1,n)$ is $n$, therefore the orbifold 
is Gorenstein.
 
\medskip
The paper is organized as follows. In Section 2 we recall some general 
facts about weighted projective spaces.  In Section 3 we 
state the main result (Thm. \ref{th2})
and we write a recipe 
which we follow for the proof. 
We prove Theorem 
\ref{th2} in Section \ref{sec:P1344} using the computation
of the Gromov-Witten invariants done in Section \ref{sec:GWI}.  
To compute those invariants we use the theory of the deformations 
of surfaces with rational double points and the deformation invariance
property of the Gromov-Witten invariants.
In Section \ref{sec:p11n} we prove that, 
for $\IP(1, \ldots ,1,n)$, the  Chen-Ruan cohomology ring is isomorphic to 
the cohomology ring of its crepant resolution (Proposition  \ref{p11n}).

\subsection*{Acknowledgments} 
Part of the work was done during the visit of two of the authors
(S.B. and F.P.) at SISSA in Trieste and hospitality and support
are gratefully acknowledged. Particular thanks go to B. Fantechi 
and Y. Ruan for very useful discussions. 
 
 \section{Weighted projective spaces} 
 \label{sec:weight-proj-spac} 
In this section we recall some basic facts about weighted projective spaces.  
 
Let $n\geq 1$ be an integer and $w=(w_0,\ldots,w_{n})$ a sequence of 
 integers greater or equal than one. Consider the action of the 
 multiplicative group $\IC^{\star}$ on $\IC^{n+1}-\{0\}$ given by: 
 $$
 \lambda\cdot(x_{0}, \ldots,x_{n}):=(\lambda^{w_{0}}x_{0}, \ldots 
 ,\lambda^{w_{n}}x_{n}). 
 $$
 The \emph{weighted projective space} $\IP(w)$ is defined as the 
 quotient stack $[\IC^{n+1}-\{0\}/\IC^{\star}]$. It is a smooth 
 Deligne-Mumford stack whose coarse moduli space, denoted $|\IP(w)|$, 
 is a projective variety of dimension $n$.  
 
According to \cite{BCS}, $\IP(w)$ is a toric stack associated to the following stacky 
 fan:
\begin{equation}\label{stackyfan}
N:=\IZ^{n+1}/\sum_{i=0}^n w_{i}v_{i},\quad \beta:\IZ^{n+1}\to N, \quad \Si,
 \end{equation}
where $v_0,...,v_n$ is the standard basis of $\IZ^{n+1}$, $\be$ is the  canonical projection, 
and $\Si \subset N\ot_{\IZ}\IQ$ is the fan whose cones
are generated by any proper subset of \\
$\{ \be(v_0)\ot 1,...,\be(v_n)\ot 1 \}$.  
 
\medskip
 The weighted projective space $\IP(w)$ comes with a natural invertible 
 sheaf  $\cO_{\IP(w)}(1)$ defined as follows: for any scheme $Y$ and any morphism 
 $Y\to \IP(w)$ given by a principal $\IC^{\star}$-bundle $P\to Y$ and a 
 $\IC^{\star}$-equivariant morphism $P\to \IC^{n+1}-\{0\}$,  
 $\cO_{\IP(w)}(1)_{Y}$ is the sheaf of sections of the associated line bundle of 
 $P$. 

\bigskip 
We will use the following 
\begin{notation}\nf
An \textit{orbifold} is a smooth algebraic Deligne-Mumford stack
over $\IC$ with generically trivial stabilizers.
The orbifold is said to be \textit{Gorenstein}
if all its ages are integers.
\end{notation}
 
 \begin{proposition}\label{prop:RG}\text{}  
\noindent (1) The Deligne-Mumford stack $\IP(w)$ is an orbifold
if and only if 
the greatest common divisor of $w_{0}, \ldots ,w_{n}$ is $1$.  \\
\noindent (2) The orbifold  $\IP(w)$ is Gorenstein if and only if 
 $w_{i}$ divides  $\sum_{j=0}^n w_j$ for any      $i$.  
 \end{proposition}  
  
 \begin{proof}\text{}  
(1) In general, the generic stabilizer of a toric Deligne-Mumford stack 
which is associated to the stacky fan 
     $(N,\beta,\Sigma)$ is isomorphic to the torsion part of $N$. 
This implies the first part. 

\vspace{0.2cm}

\noindent (2) For any       $i\in\{0, \ldots ,n\}$, set   
$$
U_{i}:=\{(x_{0}, \ldots , x_{n})\in\IC^{n+1}-\{0\}\mid x_{i}=1\}.
$$  
The trivial  $\IC^{\star}$-bundle on $U_{i}$ and the  
$\IC^{\star}$-equivariant morphism $U_{i}\times \IC^{\star}\to\IC^{n+1}-\{0\}$,
$ (x,\lambda)\mapsto \lambda\cdot x  $,   
define an \'etale morphism $\varphi_i :U_{i}\to\IP(w)$ such that 
$\sqcup_i \varphi_i : \sqcup_i U_i \ra \IP(w)$ is a covering.  
The group $\mcal{U}_{w_i}\subset \IC^\star$ of $w_i$-th 
roots of the unity acts linearly and diagonally on $U_i$ with weights
$w_0,...,w_n$. \\
Let now $x\in U_i$, and let ${\rm exp}\left( \fr{2\pi \ii k}{w_i}\right) \in \mcal{U}_{w_i}$
be an element that fixes $x$. Then
\begin{equation}\label{age}
{\rm age}\left( x, {\rm exp}\left( \fr{2\pi \ii k}{w_i}\right) \right) = 
\sum_{j= 0}^n \left\{ \fr{w_j k}{w_i} \right\},
\end{equation}
where the brackets $\left\{ \right\}$ means the fractional part of a rational number.
The claim follows. 
 \end{proof}

 In dimension $1$, the only weighted projective space which is 
 Gorenstein is $\IP(1,1)\cong \IP^{1}$. In dimension $2$ and $3$,  
the complete list of Gorenstein weighted 
 projective spaces is given by the following weights : 
 \begin{equation}  
 \begin{array}{|c||c c c c|}  
 \hline \mbox{Dimension } 2 & & \mbox{Dimension } 3 & & \\\hline  
 (1,1,1) & (1,1,1,1) & (1,2,2,5)& (2,3,3,4) & (2,3,10,15) \\  
 (1,1,2) & (1,1,1,3) & (1,1,4,6) & (1,2,6,9) & (1,6,14,21) \\  
 (1,2,3) & (1,1,2,2) & (1,2,3,6) & (1,4,5,10) &\\  
           & (1,3,4,4) & (1,1,2,4) & (1,3,8,12) &\\  
 \hline  
 \end{array}  
 \label{tableRG}  
 \end{equation} 
 In dimension $n$, the problem of determining all 
 Gorenstein $\IP(w)$ is equivalent  
 to the problem of \emph{Egyptian fractions}, \ie the number of solutions of  
 $1=\frac{1}{x_0}+\cdots+\frac{1}{x_n}$ with $1\leq x_0\leq\ldots\leq x_n$ (see  
 \cite{Egypt}). Hence, there is a finite number of such $\IP(w)$.

 \section{The main result 
}\label{sec:CCRC}  
In this section we state our main result, Theorem \ref{th2}.
Then we present a recipe that we will follow during
the proof. 

Recall that, given a Gorenstein orbifold $\cX$,
a resolution of singularities $\rho :Z \ra X$ is crepant if $\rho^*K_X \cong K_Z$.   

\medskip
In order to state the theorem, we recall the definition 
of the quantum corrected cohomology ring \cite{Ru}.

Let $\cX$ be a Gorenstein orbifold with projective
coarse moduli space $X$. 
Let  $\rho :Z \ra X$ be a crepant resolution such that 
$Z$ is projective. Let $\rm{N}^+(Z)\subset \rm{A}_1(Z;\IZ)$ be the 
 monoid of effective $1$-cycles in $Z$, and set 
$$ 
 \rm{M}_{\rho}(Z):= {\rm Ker }(\rho_{\star}) \cap  N^+(Z), 
 $$ 
where $\rho_{\star}:{\rm A}_{\star}(Z;\IZ) \ra {\rm A}_{\star}(X;\IZ)$
is the morphism of Chow groups induced by the map $\rho$.

We assume that $\rm{M}_{\rho}(Z)$ is generated by a finite number of classes
of rational curves which are linearly independent over $\IQ$.
We fix a set of such generators: $\Gamma_1,...,\Gamma_m$.
Then, any  $\Ga\in \rm{M}_{\rho}(Z)$ can be written in a unique way as:  
 $$\Ga=\sum_{\ell=1}^m d_\ell \Gamma_\ell,$$  
for some  non negative integers $d_\ell$. 

\medskip 
The \textit{quantum corrected cohomology ring} of $Z$
is defined as follows.
We assign a formal variable $q_\ell$ for each $\Gamma_\ell$, hence 
$\Ga = \sum_{\ell=1}^m d_\ell \Gamma_\ell\in \rm{M}_{\rho}(Z)$ corresponds to the monomial $q_1^{d_1}\cd \cd \cd q_m^{d_m}$.  
 The \textit{quantum $3$-points function} is 
by definition:  
         \begin{equation}\label{qc3}  
         \lan \al_1, \al_2, \al_3 \ran_{\rm q}(q_1,...,q_m):= \sum_{d_1,...,d_m > 0}  
         \Psi_{\Ga}^Z(\al_1, \al_2, \al_3 )q_1^{d_1}\cd \cd \cd q_m^{d_m},  
         \end{equation}  
where $\al_1,\al_2,\al_3\in H^\star (Z;\IC)$ and 
$\Psi_{\Ga}^Z(\al_1, \al_2, \al_3 )$ is the Gromov-Witten invariant of $Z$
of genus zero, homology class $\Ga$ and three marked points.  
  
\medskip
 We assume that \eqref{qc3} defines an analytic function of the variables  
 $q_1,...,q_m$ on some region of the complex space $\IC^m$, it will be denoted  
 by $\lan \al_1, \al_2, \al_3 \ran_{\rm q}$. In the following, when  we  evaluate  
 $\lan \al_1, \al_2, \al_3 \ran_{\rm q}$ on a point $(q_1,...,q_m)$,  we will  
 implicitly assume that it is defined on such a point.  

\medskip
The \textit{quantum corrected triple intersection}  is defined by:  
         \begin{align*}  
         \lan \al_1, \al_2, \al_3 \ran_{\rho}(q_1,...,q_m):  
         = \int_Z \al_1 \cup \al_2 \cup \al_3 +\lan \al_1, \al_2, \al_3 \ran_{\rm q}
(q_1,...,q_m).
         \end{align*}  
The   \textit{quantum corrected cup product}  $\al_1 \ast_{\rho} \al_2$ 
of two classes $\al_1, \al_2 \in H^{*}(Z;\IC)$ is defined by  
 requiring that:  
         \begin{align*}  
         \int_Z (\al_1 \ast_{\rho} \al_2) \cup \al =
\lan \al_1 ,\al_2, \al \ran_{\rho}(q_1,...,q_m)\quad  
          \quad \forall\al \in H^\star (Z;\IC).  
         \end{align*}  
  
 The following result holds.  
 \begin{proposition}\label{qccohomology}  
 For any $(q_1,...,q_m)$ belonging to the domain of the quantum  
 $3$-points function, the quantum corrected cup product $\ast_{\rho}$ satisfies  
 the following properties:  
         \begin{description}  
         \item[\nf{Associativity}] it is associative on $H^\star(Z;\IC)$,
           moreover it has a  
           unit which coincides with the unit of the usual cup product
           of $Z$.  
         \item[\nf{Skewsymmetry}] 
           $\al_1 \ast_{\rho} \al_2 = (-1)^{{\rm deg}\al_1 \cd{\rm deg}\al_2} \al_2\ast_{\rho}\al_1$  
                 for any $\al_1,\al_2 \in H^\star(Z;\IC)$.  
         \item[\nf{Homogeneity}] for any $\al_1,\al_2 \in H^\star(Z;\IC)$,  
           $\rm{deg}~(\al_1\ast_{\rho} \al_2 )= 
           \rm{deg}~\al_1 + \rm{deg}~\al_2$.  
         \end{description}  
 \end{proposition}  
  
\medskip
For any $(q_1,...,q_m)$ belonging to the domain of the quantum corrected  
$3$-point function, the resulting ring $\left( H^\star (Z;\IC), \ast_{\rho} \right)$
is the \textit{quantum corrected cohomology ring}
with the quantum parameters specialized at $(q_1,...,q_m)$,
and it will also be denoted as $H^\star_{\rho}(Z;\IC)(q_1,...,q_m)$.

\medskip
Let now $\cX :=\IP(1,3,4,4)$, let $\rho:Z\ra X$
be the crepant resolution defined in Sec.\ref{sec:P1344} (1)
and let $\Ga_1,\Ga_2,\Ga_3,\Ga_4 \in {\rm M}_\rho(Z)$
be defined in Sec. \ref{sec:P1344} (3). \\
Then we have the following
\begin{theorem} \label{th2} 
  For  $(q_1,q_2,q_3,q_4) \in \{ (\ii,\ii,\ii,0), (-\ii,-\ii,-\ii,0) \}$ 
  there is a ring isomorphism 
\begin{equation*} 
H^{\star}_{\rho}(Z;\IC)(q_1,q_2,q_3,q_4) \cong H^{\star}_{\rm CR}(\IP(1,3,4,4);\IC) 
\end{equation*} 
which is an isometry with respect to the Poincar\'e pairing
on $H^\star_{\rho}(Z;\IC)(q_1,q_2,q_3,q_4)$ 
and the Chen-Ruan pairing on $H^\star_{\rm CR}(\IP(1,3,4,4);\IC)$.\\
Furthermore, in each case an explicit
isomorphism is given by the linear map \eqref{ri} and \eqref{ri2}
respectively.
\end{theorem} 

\bigskip
\begin{remark}\nf
In general, the \textit{Cohomological Crepant Resolution
Conjecture} asserts that there exist roots of the unity
$c_1,...,c_m$ and a ring isomorphism (see \cite{Ru})
 $$  
 H^\star_{\rho}(Z;\IC)(c_1,...,c_m) \cong H^\star_{\rm CR}({\cX};\IC).  
 $$  

When the orbifold $\cX$ satisfies the hard Lefschetz condition, 
Ruan's conjecture has been generalized by Bryan and Graber \cite{BG}:
the  Frobenius manifolds associated to the big quantum cohomology of $Z$
and of $\cX$ are analytic continuations of each other. Bryan and Graber's conjecture
is called the  \textit{Crepant Resolution Conjecture}, see 
\cite{BG}, \cite{BGP}, \cite{CCIT}, \cite{CCIT2}, \cite{Wise} for 
more details and its verification in some examples.

Surprisingly, we obtain that some quantum parameters 
can be put to zero. This result is strange in regard
to the Conjecture. One can observe that in our
computations for $\IP(1,3,4,4)$ the quantum parameters corresponding to 
the transversal ${\rm A}_3$ singularity 
are evaluated at primitive $4$-th roots of the unity,
as predicted by the conjecture,
and the quantum parameter evaluated
to zero corresponds to the singular point
$\fr{1}{3}(1,1,1)$. Note that, also in the case of
$\IP(1,...,1,n)$ we set the quantum parameter
to zero (see Prop. \ref{p11n}) and the singularity
is of type $\fr{1}{n}(1,...,1)$.
\end{remark}
 
\medskip
We present below the recipe that 
we will follow. We write 
it for a general Gorenstein orbifold $\IP(w)$
because it could be used in a more general case. 

\begin{notation}\nf
For what concerns Chow rings, homology and cohomology of
Deligne-Mumford stacks, we refer the reader to \cite{AGVgwdms} Sec. 2;
we follow the same notations. 
\end{notation}

We begin with a general result.
\begin{lemma}\label{cono relativo}
Let $\IP(w)$ be a Gorenstein orbifold
and let  $\rho :Z \ra | \IP(w) |$ 
be a crepant resolution associated to a subdivision $\Si'$ of $\Si$
and the identity morphism of $N$.
Then the cone $\rm{M}_{\rho}(Z)$ is polyhedral.
\end{lemma}
\begin{proof}
Let $\Si'(n-1)$ be the set of $(n-1)$-dimensional cones of $\Si'$. Then
$$
\rm{M}_{\rho}(Z) = \left\{ \sum_{\nu \in \Si'(n-1)}\ga_{\nu}\left[{\rm V}(\nu)\right] \, \mid \,
 \ga_{\nu} \in \IN, \, 
\rho_{\star}\left(\sum_{\nu \in \Si'(n-1)}\ga_{\nu}\left[{\rm V}(\nu)\right]\right)=0\right\},
$$
where, for any $\nu \in \Si'(n-1)$, ${\rm V}(\nu)$
denotes the rational curve in $Z$ stable under the 
torus action which is associated  to $\nu$ \cite{F}, $\left[{\rm V}(\nu)\right]$
is the induced Chow class.\\
Let now $L \in {\rm Pic}(|\IP(w)|)$ be an ample line bundle.
From standard intersection theory we have (see e.g. \cite{Fit}):
$$
\rho_{\star}\left(\sum_{\nu \in \Si'(n-1)}\ga_{\nu}\left[{\rm V}(\nu)\right]\right)=0
\quad \mb{if and only if} \quad 
c_1\left(\rho^\star L\right)\cap\left(\sum_{\nu \in \Si'(n-1)}\ga_{\nu}\left[{\rm V}(\nu)\right]\right)
=0.
$$
Since $c_1\left(\rho^\star L\right)\cap \left[{\rm V}(\nu)\right]\geq 0$ for any $\nu$, 
it follows that 
$$
\rm{M}_{\rho}(Z) = \left\{\sum_{\nu \in \Si'(n-1)}\ga_{\nu}\left[{\rm V}(\nu)\right] \, \mid \,
 \ga_\nu \in \IN, \, \rho_{\star}\left(\left[{\rm V}(\nu)\right]\right)=0\right\},
$$
hence the claim.
\end{proof}

\bigskip 
The steps that we will follow during the proof of our
results are the following.
 \begin{enumerate}  
 \item Let $(N,\be, \Si)$ be the stacky fan \eqref{stackyfan},
choose a subdivision $\Si'$
of $\Si$ such that the  morphism $\rho : Z \ra |\IP(w)|$
associated to $\Si'$ and the identity of $\rm N$
is a crepant resolution, where $Z$ is the toric variety
associated to $\Si'$.\\
In dimension $2$, this is the classical  
Hirzebruch-Jung algorithm. In dimension $3$, this is always possible (see \cite{CR}). 
\medskip
\item Set 
$$
H:= c_1\left( \cO_{\IP(w)}(1) \right) \in H^2(\IP(w);\IC ),
$$
and set 
$$
h:=\rho^\star H \in H^2(Z;\IC).
$$
Note that we identify $ H^2(\IP(w);\IC )$ with $ H^2(|\IP(w)|;\IC )$.
 Let us denote by
$b_i \in H^2(Z;\IC)$ (by $e_1,...,e_d\in H^2(Z;\IC)$ resp.)
the first Chern class of the line bundle associated to the torus invariant
divisor corresponding to the ray of $\Si'$ generated by $\be(v_i)$ 
(the rays in $\Si'(1)-\Si(1)$ resp.), for any $i\in \{0,...,n\}$.\\
Then compute  a presentation of $H^\star(Z;\IC)$ 
as a quotient of  
 $\IC[h,e_1,\ldots,e_d]$, e.g. following \cite{F} Sec. 5.2.
 
Note that, since $\rho$ is crepant, we have the following equality:
\begin{equation}\label{h}
h=\fr{1}{\sum_{i=0}^n w_i}\left(\sum_{i=0}^nb_i + \sum_{j=1}^de_j\right).
\end{equation}

\medskip
\item Find the generators of ${\rm M}_\rho(Z)$. 
Here is a possible way:
$H$ is an ample line bundle (see e.g. \cite{F} Sec. 3.4),
then, following the proof of Lemma \ref{cono relativo},
we see that  $ {\rm M}_\rho(Z)$
is generated by the set 
\begin{equation}\label{generatori cono relativo}
\left\{ [V(\nu)] \, \mid \, h\cap \left[V(\nu)\right]=0, \, \nu \in \Sigma'(n-1)-\Sigma(n-1)\right\}.
\end{equation}
We remark that, given the fans $\Si$ and $\Si'$ and using \eqref{h}, 
the determination of the set \eqref{generatori cono relativo} is straightforward
using for instance the results in \cite{F} Sec. 5.2.

\medskip
 \item Compute the Gromov-Witten invariants in a basis of $H^\star(Z;\IC)$
and express the quantum corrected cohomology ring $H^\star_\rho(Z;\IC)(q_1,\ldots,q_m)$
in such a basis.  
\medskip
   \item Compute a presentation of $H^\star_{\CR}(\IP(w);\IC)$.
Note that using the combinatorial model presented in \cite{BMPmodel}, or results from \cite{M},
one obtains directly  
a basis and the multiplicative table of $H^\star_{\CR}(\IP(w);\IC)$. 
Equivalently the results from \cite{BCS} and \cite{CCLT} 
can be used to obtain a presentation
of the Chen-Ruan cohomology ring of weighted projective spaces.

\medskip
   \item Find a suitable evaluation of the quantum
parameters and a linear map
$H^\star_\rho(Z;\IC)(c_1,\ldots,c_m)\ra H^\star_{\CR}(\IP(w);\IC)$
that induces a  ring isomorphism.  
 \end{enumerate}  

\medskip
\begin{example}\nf
As an example we work out the steps above for $\cX:=\IP(1,1,2,2)$.
Note that $\IP(1,1,2,2)$ is an orbifold with transversal ${\rm A}_1$
singularities, hence Ruan's conjecture in this case follows from \cite{P}.

\medskip
(1) Identify the stacky fan $(N,\be,\Si)$ with 
$(\IZ^3,\{\Lambda(\be(v_i))\}_{i\in \{0,1,2,3\}}, \Si)$, where
the $v_i$ are defined in \eqref{stackyfan} and 
$\Lam :N\ra \IZ^3$ is the isomorphism defined by sending
$v_0 \mapsto (-1,-2,-2)$, $v_1 \mapsto (1,0,0)$, $v_2 \mapsto (0,1,0)$ 
and $v_3 \mapsto (0,0,1)$.

The subdivision $\Si'$ of $\Si$ is obtained  by adding the ray generated by
$$
{\rm P}:=(0,-1,-1)=\frac{1}{2}\left(\Lam(\be(v_0))+\Lam(\be(v_1))\right)
$$
as shown in Figure \ref{fig3}.  
\begin{figure}[!ht]  
\begin{center}  
\psfrag{A}{$\Lam(\be(v_{1}))$} \psfrag{B}{$\Lam(\be(v_{2}))$} 
\psfrag{C}{$\Lam(\be(v_{3}))$} \psfrag{D}{$\Lam(\be(v_{0}))$}  
\psfrag{E}{${\rm P}$}  \includegraphics[width=0.4\linewidth]{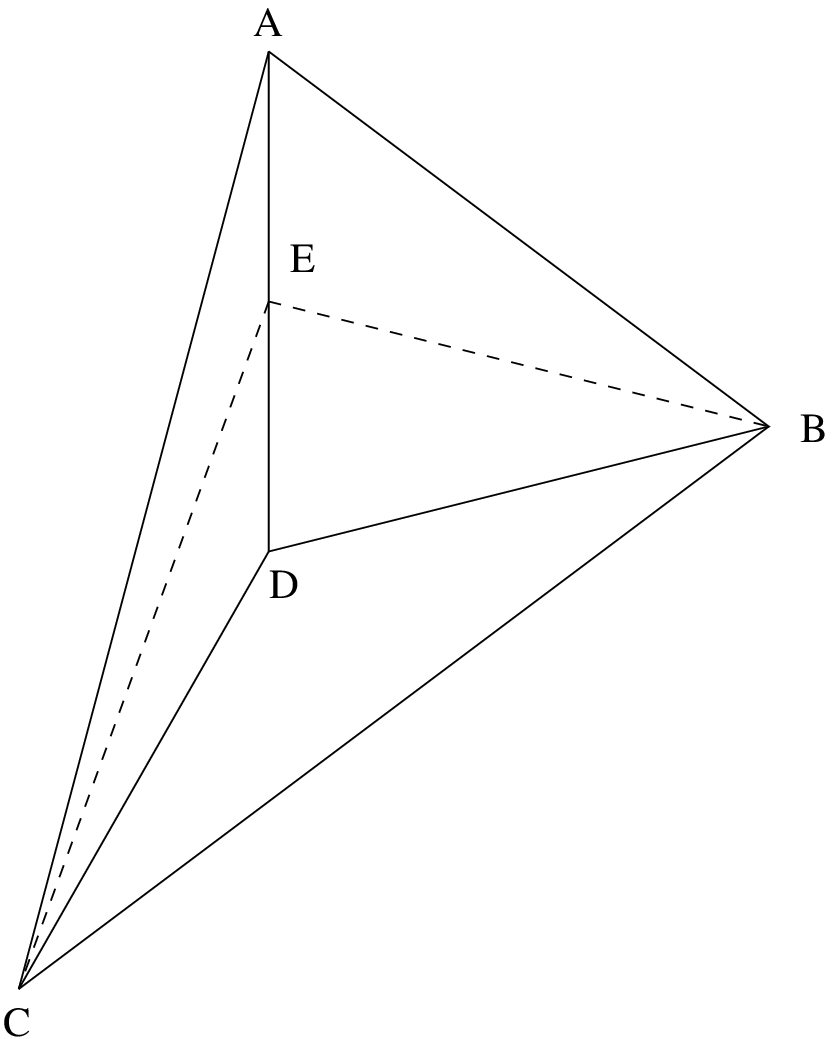}  
\end{center}\caption{Polytope of $\IP(1,1,2,2)$ and a crepant resolution.}\label{fig3}  
\end{figure} 

\medskip 
(2) We have $ H^\star(Z;\IC)\cong \IC[h,e]/\langle h^2+\frac{1}{4}e^2-he,h^2e\rangle$.  

\medskip
(3) The new effective curves are $\PD(eb_0)$ and $\PD(eb_2)$ but only  
$\PD(eb_2)$ is contracted by $\rho$, so ${\rm M}_\rho(Z)$ is generated by 
$\Ga_1:=\PD(eb_2)$.   Note that $\PD(eb_2)=\PD(2he)$.

\medskip
(4) Since  $\dim \left( [\overline{\cM}_{0,3}(Z,d\Gamma_1)]^{\rm vir}\right) =3$, the Gromov-Witten invariants  
$\Psi^Z_{d\Gamma_1}(\al_1,\al_2,\al_3)$ are not zero only when all the $\al_i$  
have degree $2$. By the Divisor Axiom we have:  
$$  
\Psi^Z_{d\Gamma_1}(\al_1,\al_2,\al_3)=\left(\int_{\Gamma_1}\al_1\right)  
\left(\int_{\Gamma_1}\al_2\right)\left(\int_{\Gamma_1}\al_3\right) \cdot d^3\cdot  
\Psi^Z_{\Gamma_1}(\cdot).  
$$  
From \cite[Proposition IV.3.13]{M} we get:  
$$  
\int_Zh^3=\int_Z\rho^\star c_1(\cO_{\IP(w)}(1))^3=\int^{\CR}_{\IP(1,1,2,2)}c_1(\cO_{\IP(w)}(1))^3=\frac{1}{4}.  
$$  
This gives:  
\begin{align*}  
\Psi^Z_{d\Gamma_1}(h,h,h)&=0 & \Psi^Z_{d\Gamma_1}(h,h,e)&=0\\  
\Psi^Z_{d\Gamma_1}(h,e,e)&=0 & \Psi^Z_{d\Gamma_1}(e,e,e)&=-2^3d^3\cdot \Psi^Z_{\Gamma_1}(\cdot)  
\end{align*}  
From \cite{P} we get:  
$$  
\Psi^Z_{d\Gamma_1}(\cdot)=\frac{2}{d^3}\int_{|\IP(2,2)|}c_1(\cO_{\IP(2,2)}(1))=\frac{1}{d^3}.  
$$  
This gives the quantum correction:  
$$  
e\ast_\rho e=-4h^2+\left(4+\frac{8q_1}{1-q_1}\right)he,  
$$  
hence:  
\begin{equation}\label{qccrp1122}  
H^\star_\rho(Z;\IC)(q_1)\cong \IC((q_1))[h,e]/\langle  
h^2e,h^2+\frac{1}{4}e^2-he-\frac{2q_1}{1-q_1}he\rangle.  
\end{equation}

\medskip
(5) The Chen-Ruan cohomology ring has the following presentation
(see \cite{BMPmodel} or Sec. \ref{sec:P1344} for more details):  
$$  
H^\star_{\CR}(\IP(1,1,2,2);\IC)\cong \IC[H,E]/\langle H^2-E^2,H^2E\rangle.  
$$  

\medskip
(6) Set $q_1=-1$ in \eqref{qccrp1122} and define the following map:  
\begin{eqnarray}\label{ccrcp1122}  
H^\star_{\CR}(\IP(1,1,2,2))&\ra& H_\rho^\star(Z)(-1)  \\
H&\mapsto& h,\nonumber \\  
 E&\mapsto & \frac{\ii}{2} e.\nonumber
\end{eqnarray}  
Then \eqref{ccrcp1122} is a ring isomorphism.
\end{example}

\section{Proof of Theorem \ref{th2}}\label{sec:P1344}\text{}  
In this section we prove Theorem \ref{th2} 
by following the steps above for $\cX = \IP(1,3,4,4)$.

We remark that this result confirms 
for the transversal ${\rm A}_3$ singularity case Conjecture 1.9 in \cite{P} regarding 
the values of the $q_i$'s, secondly, the change of variables is 
inspired from those of \cite{NW} (see also \cite{BGP} and \cite{CCIT2}). 
 
The coarse moduli space of $\IP(1,3,4,4)$ has a transversal $A_3$ singularity  
on the line $[0:0:x_2:x_3]$ and an isolated singularity  of type
$\frac{1}{3}(1,1,1)$ at the point $[0:1:0:0]$.  
  
\medskip
(1) We identify the stacky fan $(N,\be,\Si)$ with 
$(\IZ^3,\{\Lambda(\be(v_i))\}_{i\in \{0,1,2,3\}}, \Si)$, where
the $v_i$ are defined in \eqref{stackyfan} and 
$\Lam :N\ra \IZ^3$ is the isomorphism defined by sending
$v_0 \mapsto (-3,-4,-4)$, $v_1 \mapsto (1,0,0)$, $v_2 \mapsto (0,1,0)$ 
and $v_3 \mapsto (0,0,1)$.

\begin{figure}[!ht]  
\begin{center}  
\psfrag{A}{$\Lam(\be(v_{1}))$} \psfrag{B}{$\Lam(\be(v_{2}))$} 
\psfrag{C}{$\Lam(\be(v_{3}))$} \psfrag{D}{$\Lam(\be(v_{0}))$}  
\psfrag{E}{${\rm P}_1$} \psfrag{F}{${\rm P}_2$} \psfrag{G}{${\rm P}_3$} \psfrag{H}{${\rm P}_4$}  
 \includegraphics[width=0.4\linewidth]{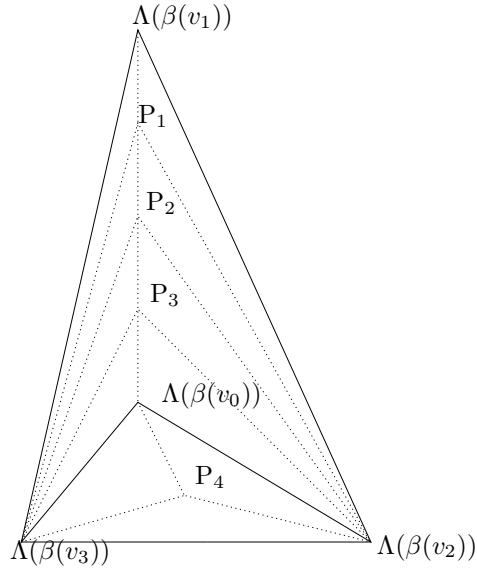}  
\end{center}\caption{Polytope of $\IP(1,3,4,4)$ and a crepant resolution}\label{fig4}  
\end{figure}

A crepant resolution of $|\IP(1,3,4,4)|$ can be constructed using standard methods
in toric geometry. More precisely, consider 
the integral points
\begin{align*}  
  {\rm P}_1&:=(0,-1,-1)=\frac{3}{4}\Lam(\be(v_{1}))+\frac{1}{4}\Lam(\be(v_{0})),\\  
{\rm P}_{2}&:=(-1,-2,-2)=\frac{1}{2}\Lam(\be(v_{1}))+\frac{1}{2}\Lam(\be(v_{0})), \\  
 {\rm P}_{3}&:=(-2,-3,-3)=\frac{1}{4}\Lam(\be(v_{1}))+\frac{3}{4}\Lam(\be(v_{0})),\\  
\mb{and} \qquad  {\rm P}_{4}&:=(-1,-1,-1)=\frac{1}{3}\Lam(\be(v_{0}))+\frac{1}{3}\Lam(\be(v_{2}))+
\frac{1}{3}\Lam(\be(v_{3})),  
\end{align*} 
then subdivide $\Si$ inserting the rays generated by ${\rm P}_1$,
${\rm P}_2$, ${\rm P}_3$ and  ${\rm P}_4$
as shown in  Figure \ref{fig4}. Let
$\Si'$ be the fan obtained after this subdivision, let $Z$ be the toric variety
and $\rho :Z \ra |\IP(1,3,4,4)|$ be the morphism associated to the identity
on $\IZ^3$. Then $Z$ is smooth and $\rho$ is crepant.
The crepancy of $\rho$ follows from the existence
of a continuous piecewise linear function $|\Si | \ra \IR$
which is linear when restricted to each cone of $\Si$
and associates the value $-1$ to the minimal lattice points of the 
rays of $\Si'$ (see \cite{F}, Sec. 3.4). 
  
\medskip
(2) We denote by 
$b_i\in H^2(Z;\IC)$ ($e_j$ resp.) the first Chern class of the line bundle 
associated to the torus invariant divisor corresponding to the ray
generated by $\Lam(\be(v_i))$ (${\rm P}_j$ resp.), for $i\in\{0,1,2,3\}$
($j\in\{1,2,3,4\}$ resp.), furthermore we define
$$
h=\fr{1}{12}\left(\sum_{i=0}^3 b_i + \sum_{j=1}^4 e_j\right).
$$
The cohomology ring of $Z$ has the following
presentation:
it is isomorphic to the quotient of the polynomial ring
$\IC[h,e_{1},e_{2},e_{3},e_{4}]$ by the ideal generated by   
\begin{align*}  
3he_{4}, \   e_{1}e_{3}, \ e_{1}e_{4}, \ e_{2}e_{4}, \ e_{3}e_{4}, \\  
e_{1}^{2} - 10 he_{1}-4he_{2}-2he_{3}+24h^{2},\\  
e_{1}e_{2} + 3 he_{1}+2 he_{2}+he_{3}-12h^{2},\\  
e_{2}^{2}-6he_{1}-12he_{2}-2he_{3}+24h^{2},\\  
e_{2}e_{3}+3he_{1}+6he_{2}+he_{3}-12h^{2},\\  
e_{3}^{2}-6he_{1}-12he_{2}-14he_{3}+24h^{2},\\  
16h^{2}e_{1},\ 16h^{2}e_{2}, \ 16h^{2}e_{3}, 16h^{3}-\frac{1}{27}e_{4}^{3}.  
\end{align*}  
We fix the following basis of the vector space
$H^\star(Z;\IC)$:\\
$$1,\, h,\, e_1,\, e_2,\, e_3,\, e_4,\,  h^2,\, he_1,\, he_2,\, he_3,\, e_4^2,\, h^3.$$
  
\medskip
(3) As explained in Sec. \ref{sec:CCRC}, the cone ${\rm M}_\rho(Z)$ 
can be directly determined from the combinatorial data $\Si$ and $\Si'$.
In our case ${\rm M}_\rho(Z)$ is generated by
$\Gamma_1:=\rm{PD}(4he_{1})$, $\Gamma_2:=\rm{PD}(4he_{2})$,  
$\Gamma_3:=\rm{PD}(4he_{3})$ and $\Gamma_4:=\rm{PD}(-\frac{1}{3}e_{4}^{2})$,
where ${\rm PD }$ means Poincar\'e dual.

\medskip  
(4) Here we give a presentation of the quantum corrected cohomology
ring $H^\star_{\rho}(Z;\IC)(q_1,q_2,q_3,0)$.

We first notice that any curve of homology class  
$d_4\Gamma_4$ is disjoint from any other curve of class $d_1\Gamma_1 +d_2 \Gamma_2 + d_3\Gamma_3$;
in other words, $\overline{\cM}_{0,0}(Z,\Ga)$ is empty
if $\Ga=\sum_{\ell=1}^4 d_{\ell}\Gamma_{\ell}$ with $d_4\cd(d_1+d_2+d_3)\not=0$.
From the degree axiom it follows that we need to consider only Gromov-Witten
invariants $\Psi_{\Ga}^{Z}(\al_1,\al_2,\al_3)$ with $\al_i \in H^2(Z;\IC)$, 
$i\in \{ 1,2,3 \}$.
Finally, applying the divisor axiom we deduce the following expression for  
the quantum $3$-point function:
\begin{eqnarray*} 
& & \lan \al_1, \al_2, \al_3 \ran_{\rm q}(q_1,q_2,q_3,q_4)\\ 
&=& \sum_{d_1,d_2,d_3>0}\left( \prod_{i=1}^3 \int_{\sum_{\ell=1}^{3}d_\ell\Gamma_\ell}\al_i\right) 
{\rm deg}\left[\overline{\cM}_{0,0}(Z,\sum_{\ell =1}^3d_\ell\Gamma_\ell)\right]^{\rm vir}q_1^{d_1}q_2^{d_2}q_3^{d_3}\\ 
&+&\sum_{d_4>0}\left( \prod_{i=1}^3 \int_{d_4\Gamma_4 }\al_i\right) 
{\rm deg}[\overline{\cM}_{0,0}(Z,d_4\Gamma_4)]^{\rm vir}q_4^{d_4}. 
\end{eqnarray*} 
Since $\int_{ \Gamma_\ell}h =0$ for any $\ell\in \{1,2,3,4\}$,
$$
h\ast_{\rho} \al = h\al \qquad \mb{for any} \quad \al \in H^\star(Z;\IC),
$$
similarly 
$$
e_i \ast_{\rho} e_4 =
\begin{cases} 
 e_i  e_4=0, \quad \mbox{if} \quad i\not=4; \\ 
 \e (q_4)e_4^2, \quad \mb{otherwise},
\end{cases} 
$$
for some function $\e (q_4)$ such that $\e (0)=1$.  

As in the isomorphism of rings that we will define later we will 
put $q_{4}=0$, we only consider classes 
$\Gamma =d_{1}\Gamma_{1}+d_{2}\Gamma_{2}+d_{3}\Gamma_{3}$ for $d_{i}\in 
\IN$. 
We set $\Gamma_{\mu\nu}:=\Gamma_{\mu}+\cdots+\Gamma_{\nu}$ for 
$\mu, \nu \in \{1,2,3\}$ and $\mu\leq \nu$.  
Using Theorem \ref{GW1344} in Section \ref{sec:GWI} we get:  
\[{\rm deg}[\overline{\cM}_{0,0}(Z,\Ga)]^{\rm vir}=   
\begin{cases}  
1/d^{3}  & \mbox{if} \quad \Gamma = d\Gamma_{\mu\nu} \quad \mbox{for}
\quad \mu \leq \nu\in\{1,2,3\}
\quad \mbox{and} \quad d \in \IN\\
0 &  \mbox{otherwise}.  
\end{cases}\]  
Hence the remaining part of the multiplicative table of $H^\star_{\rho}(Z;\IC)(q_1,q_2,q_3,0)$ 
is as follows:  
\begin{align*}  
e_{1}\ast_{\rho}e_{1}&=-24h^{2}   
+\left(10  +16\frac{q_{1}}{1-q_{1}}+4\frac{q_{1}q_{2}}{1-q_{1}q_{2}}+4\frac{q_{1}q_{2}q_{3}}{1-q_{1}q_{2}q_{3}}\right)he_{1}\\  
&\hspace{1.45cm} +\left(4+4\frac{q_{2}}{1-q_{2}}+4\frac{q_{1}q_{2}}{1-q_{1}q_{2}}+4\frac{q_{2}q_{3}}{1-q_{2}q_{3}}  
 +4\frac{q_{1}q_{2}q_{3}}{1-q_{1}q_{2}q_{3}}\right)he_{2}\\  
& \hspace{1.45cm}+  
\left(2+4\frac{q_{2}q_{3}}{1-q_{2}q_{3}}  
  +4\frac{q_{1}q_{2}q_{3}}{1-q_{1}q_{2}q_{3}}\right)he_{3},  
\end{align*}  
\begin{align*}  
e_{1}\ast_{\rho}e_{2}&=12h^{2}   
+\left(-3 -8\frac{q_{1}}{1-q_{1}} +4\frac{q_{1}q_{2}}{1-q_{1}q_{2}}\right)he_{1}\\  
& \hspace{1.45cm} +\left(-2-8\frac{q_{2}}{1-q_{2}}+4\frac{q_{1}q_{2}}{1-q_{1}q_{2}}-4\frac{q_{2}q_{3}}{1-q_{2}q_{3}}\right)he_{2}\\  
&\hspace{1.45cm}+  
\left(-1-4\frac{q_{2}q_{3}}{1-q_{2}q_{3}}\right)he_{3},  
\end{align*}  
\begin{align*}  
e_{1}\ast_{\rho}e_{3}&=\left(-4\frac{q_{1}q_{2}}{1-q_{1}q_{2}}+4\frac{q_{1}q_{2}q_{3}}{1-q_{1}q_{2}q_{3}}\right)he_{1}\\  
&\hspace{0.45cm}  
+\left(4\frac{q_{2}}{1-q_{2}}-4\frac{q_{1}q_{2}}{1-q_{1}q_{2}}-4\frac{q_{2}q_{3}}{1-q_{2}q_{3}}+  
4\frac{q_{1}q_{2}q_{3}}{1-q_{1}q_{2}q_{3}}\right)he_{2}\\  
& \hspace{0.45cm}+  
\left(-4\frac{q_{2}q_{3}}{1-q_{2}q_{3}}+4\frac{q_{1}q_{2}q_{3}}{1-q_{1}q_{2}q_{3}}\right)he_{3},  
\end{align*}  
\begin{align*}  
e_{2}\ast_{\rho}e_{2}&=-24h^{2}  +\left(6  
  +4\frac{q_{1}}{1-q_{1}}+4\frac{q_{1}q_{2}}{1-q_{1}q_{2}}\right)he_{1}\\  
& \hspace{1.45cm}  
+\left(12+16\frac{q_{2}}{1-q_{2}}+4\frac{q_{1}q_{2}}{1-q_{1}q_{2}}+4\frac{q_{2}q_{3}}{1-q_{2}q_{3}}\right)he_{2}\\  
& \hspace{1.45cm} +  
\left(2+4\frac{q_{3}}{1-q_{3}}+4\frac{q_{2}q_{3}}{1-q_{2}q_{3}}\right)he_{3},  
\end{align*}  
\begin{align*}  
e_{2}\ast_{\rho}e_{3}&=12h^{2}    
+\left(-3 -4\frac{q_{1}q_{2}}{1-q_{1}q_{2}}\right)he_{1}\\  
& \hspace{1.18cm}  
+\left(-6-8\frac{q_{2}}{1-q_{2}}-4\frac{q_{1}q_{2}}{1-q_{1}q_{2}}+4\frac{q_{2}q_{3}}{1-q_{2}q_{3}}\right)he_{2}\\& \hspace{1.18cm}+  
\left(-1-8\frac{q_{3}}{1-q_{3}}+4\frac{q_{2}q_{3}}{1-q_{2}q_{3}}\right)he_{3},  
\end{align*}  
\begin{align*}  
e_{3}\ast_{\rho}e_{3}&=-24h^{2}  
+\left(6+4\frac{q_{1}q_{2}}{1-q_{1}q_{2}} +4\frac{q_{1}q_{2}q_{3}}{1-q_{1}q_{2}q_{3}}\right)he_{1}\\  
& \hspace{1.45cm} +\left(12+4\frac{q_{2}}{1-q_{2}}+4\frac{q_{1}q_{2}}{1-q_{1}q_{2}}+4\frac{q_{2}q_{3}}{1-q_{2}q_{3}}  
 +4\frac{q_{1}q_{2}q_{3}}{1-q_{1}q_{2}q_{3}}\right)he_{2}\\& \hspace{1.45cm}+  
\left(14+16\frac{q_{3}}{1-q_{3}}+4\frac{q_{2}q_{3}}{1-q_{2}q_{3}}+4\frac{q_{1}q_{2}q_{3}}{1-q_{1}q_{2}q_{3}}\right)he_{3}.  
\end{align*}  
  
\bigskip
(5) To compute the Chen-Ruan cohomology ring $H^{\star}_{\rm CR}(\cX;\IC)$
we follow \cite{BMPmodel}. The twisted sectors are indexed by the set
${\rm T}:= \left\{ {\rm exp}(2\pi \ii \ga)\, | \, 
\ga \in \left\{ 0, \frac{1}{3}, \frac{2}{3}, \frac{1}{4}, \frac{1}{2},
\frac{3}{4}\right\} \right\}$. For any $g \in {\rm T}$,
$\cX_{(g)}$ is a weighted projective space: set 
${\rm I}(g):=\left\{ i\in \{ 0,1,2,3 \} \, \mid \, g^{w_i} =1  \right\}$,
then $\cX_{(g)}=\IP(w_{{\rm I}(g)})$, where
$(w_{{\rm I}(g)})=(w_i)_{i\in {\rm I}(g)}$. The inertia stack is the disjoint union 
of the twisted sectors:
$$
{\rm I}\cX = \sqcup_{g\in {\rm T}}\IP(w_{{\rm I}(g)}).
$$ 
As a vector space, the Chen-Ruan cohomology is
the cohomology of the inertia stack; the graded structure is obtained
by shifting the degree of the cohomology of any twisted sector by 
twice the corresponding age (see \eqref{age}).
We have
\begin{eqnarray}\label{crcohomology}
H^p_{\rm CR}(\cX;\IC) &=&
 \op_{g\in {\rm T}}H^{p-2{\rm age}(g)}(\IP(w_{{\rm I}(g)});\IC)\\
&=&H^p(\IP(1,3,4,4);\IC) \op H^{p-2}(\IP(3);\IC) \op H^{p-4}(\IP(3);\IC)\op
\nonumber \\
&& H^{p-2}(\IP(4,4);\IC)\op H^{p-2}(\IP(4,4);\IC)\op H^{p-2}(\IP(4,4);\IC).
\nonumber
\end{eqnarray}
A basis of $H^{\star}_{\rm CR}(\cX;\IC)$ is easily obtained in the following
way: 
set 
$$
H, E_1, E_2, E_3, E_4\in H^{\star}_{\rm CR}(\cX;\IC) 
$$
be the image of 
$
c_1 \left(\cO_{\cX}(1)\right) \in H^{2}(\cX;\IC),
$
$1\in H^0(\cX_{\left({\rm exp}(\pi\ii /2)\right)};\IC)$,
$1\in H^0(\cX_{\left({\rm exp}(\pi\ii )\right)};\IC)$,
$1\in H^0(\cX_{\left({\rm exp}(\pi\ii 3/2)\right)};\IC)$
and $1\in H^0(\cX_{\left({\rm exp}(2\pi\ii /3)\right)};\IC)$ respectively,
under the inclusion $H^{\star}(\IP(w_{{\rm I}(g)}))\ra H^\star_{\rm CR}(\cX)$
determined by the decomposition \eqref{crcohomology}.
As a $\IC$-algebra, the Chen-Ruan cohomology ring is generated by 
$H,E_1,E_2,E_3,E_4$ with relations (see \cite{BMPmodel}):
\begin{align*}  
&& HE_{4}, \ E_{1}E_{1}-3HE_{2}, \ E_{1}E_{2}-3HE_{3}, \  
E_{1}E_{3}-3H^{2}, \\  &&
E_{2}E_{2}-3H^{2}, \ E_{2}E_{3}-HE_{1}, \ E_{3}E_{3}-HE_{2}, \  
16H^{3}-E_4^{3}, \\ &&
H^2E_1, \, H^2E_2, \, H^2E_3, \, E_1E_4, \, 
E_2E_4, \, E_3E_4. 
\end{align*}  
We see that the following elements form a basis of 
$H^\star_{\rm CR}(\cX;\IC)$ which we fix for the rest 
of the paper:
$$
1,\, H,\, E_1,\, E_2,\, E_3,\, E_4,\, H^2,\, HE_1,\, HE_2,\, HE_3,\, E_4^2,\, H^3.
$$

\begin{remark}\nf
Note that the elements of our basis 
are different from those used in \cite{BMPmodel}
by a combinatorial factor.

Other methods are suitable in order to compute the Chen-Ruan cup 
product of weighted projective spaces, here are a few:
the results in \cite{BCS}
provide a presentation of the Chen-Ruan cohomology ring
for toric Deligne-Mumford stack; results from \cite{M} and from \cite{CCLT}.
\end{remark}
  
\medskip
(6) We have shown in (4) that, for cohomology classes $\al_1$ and $\al_2$,
the product $\al_1\ast_{\rho}\al_2\in H^\star_\rho(Z;\IC)(q_1,q_2,q_3,0)$ 
differs from the usual cup product only if $\al_1,\al_2\in\{e_1,e_2,e_3\}$. 
We now set $q_1=q_2=q_3=\ii$ and write $e_i\ast_{\rho}e_j$ with respect to
the basis of $H^\star(Z;\IC)$ fixed in point (2),
we have:  
\begin{align*}  
&e_{1}\ast_{\rho}e_{1}=  
-24h^{2}+(-2+6\ii)he_{1}-4he_{2}+(-2-2\ii)he_{3},\\  
& e_{1}\ast_{\rho}e_{2}=12h^{2}+(-1-4\ii)he_{1}+(2-4\ii)he_{2}+he_{3},\\  
& e_{1}\ast_{\rho}e_{3}=  
-2\ii he_{1}-2\ii he_{3},\\  
& e_{2}\ast_{\rho}e_{2}=  
-24h^{2}+(2+2\ii)he_{1}+8\ii he_{2}+(-2+2\ii)he_{3},\\  
& e_{2}\ast_{\rho}e_{3}=  
12h^{2}-he_{1}+(-2-4\ii)he_{2}+(1-4\ii)he_{3},\\  
& e_{3}\ast_{\rho}e_{3}=  
-24h^{2}+(2-2\ii)he_{1}+4he_{2}+(2+6\ii)he_{3}.  
\end{align*}  
We define a linear map   
\begin{equation}\label{ri}
H_\rho^\star(Z;\IC)(\ii,\ii,\ii,0)\rightarrow H^\star_{\CR}(\IP(1,3,4,4);\IC)
\end{equation} 
as follows: we send
\begin{equation*}  
\begin{array}{cccccccc}
\left( 
\begin{array}{c}
h\\ e_1\\e_2\\e_3\\e_4
\end{array}
\right)
\begin{array}{c}
\longmapsto
\end{array}
\left(
\begin{array}{ccccc}
1 & 0 & 0 & 0 & 0 \\
0 & -\sqrt{2} & -2\ii & \sqrt{2} & 0 \\
0 & -\ii\sqrt{2} & 2\ii & -\ii\sqrt{2} & 0 \\
0 & \sqrt{2} & -2\ii & -\sqrt{2} & 0 \\
0 & 0 & 0 & 0 & 3{\rm exp}\left(\fr{2\pi\ii}{3}\right)
\end{array}
\right)
\left(
\begin{array}{c}
H\\ E_1\\ E_2\\
E_3\\ E_4
\end{array}
\right),
\end{array}  
\end{equation*} 
the image of the other elements of the basis is uniquely determined 
by requiring that \eqref{ri} is a ring isomorphism.
  
A direct computation shows that \eqref{ri} is a ring isomorphism 
and that it is an isometry with respect to the inner products given by the 
Poincar\'e duality and the Chen-Ruan pairing respectively.

\medskip
The case where $q_{1}=q_{2}=q_{3}=-\ii$ and $q_{4}=0$ is analogous to the previous one.
We define a linear map
\begin{equation}\label{ri2}
H_\rho^\star(Z;\IC)(-\ii,-\ii,-\ii,0)\rightarrow H^\star_{\CR}(\IP(1,3,4,4);\IC)
\end{equation}
by sending 
\begin{equation*}  
\begin{array}{cccccccc}
\left( 
\begin{array}{c}
h\\ e_1\\e_2\\e_3\\e_4
\end{array}
\right)
\begin{array}{c}
\longmapsto
\end{array}
\left(
\begin{array}{ccccc}
1 & 0 & 0 & 0 & 0 \\
0 & -\sqrt{2} & 2\ii & \sqrt{2} & 0 \\
0 & \ii\sqrt{2} & -2\ii & \ii\sqrt{2} & 0 \\
0 & \sqrt{2} & 2\ii & -\sqrt{2} & 0 \\
0 & 0 & 0 & 0 & 3{\rm exp}\left(\fr{2\pi\ii}{3}\right)
\end{array}
\right)
\left(
\begin{array}{c}
 H\\ E_1\\ E_2\\
E_3\\ E_4
\end{array}
\right)
\end{array}  
\end{equation*} 
and extending to the remaining part of the basis in the unique way such that
the resulting map is a ring isomorphism.

Also in this case a direct computation shows that \eqref{ri2} is a ring isomorphism
and it respects the inner pairings. 

\begin{remark}\nf
Note that the isomorphisms \eqref{ri} and \eqref{ri2} are those conjectured in \cite{P}.
\end{remark}

\section{ Gromov-Witten invariants of the resolution of 
  $|\IP(1,3,4,4)|$} \label{sec:GWI} 
 
In this section we compute the  Gromov-Witten invariants of 
the crepant resolution of $|\IP(1,3,4,4)|$ of genus $0$, homology class 
$\Gamma= d_1\Gamma_1 +d_2 \Gamma_2 +d_3\Gamma_3$
and without marked points.  Our result confirms 
Conjecture 5.1 \cite{P}. We follow the notations from Section 
\ref{sec:P1344}. 
\begin{theorem}\label{GW1344}  
Let $\rho :Z\ra |\IP(1,3,4,4)|$ be the crepant resolution of 
$|\IP(1,3,4,4)|$ defined in Section \ref{sec:P1344} (1), and let   
$\Ga= d_1\Gamma_1 +d_2 \Gamma_2 +d_3\Gamma_3$. Then   
\[ {\rm deg}[\overline{\cM}_{0,0}(Z,\Ga)]^{\rm vir}=  
\begin{cases}  
{1}/{d^3}\quad \mb{if} \quad \Ga= d\sum_{i=\mu}^{\nu}\Gamma_i, \quad \mbox{with}
\quad \mu\leq\nu \in \{1,2,3\};\\  
0 \quad \mb{otherwise}.  
\end{cases}\]  
\end{theorem}

To prove Theorem \ref{GW1344} we use the deformation invariance property of 
the Gromov-Witten invariants. More precisely: we define an open 
neighborhood $V$ of the singular locus $|\IP(4,4)| \subset |\IP(1,3,4,4)|$,
we construct an 
explicit deformation of $V$ and then a simultaneous resolution.  
This gives a deformation of $\rho^{-1}(V)$, a neighborhood  of the component of 
the exceptional divisor which lies over $|\IP(4,4)|$.
We will denote this deformation by $\overline{{\rm Graph} (\mu)}_{t}$, $t\in \De$.
Next we relate the Gromov-Witten invariants of $Z$ we are interested in 
with some Gromov-Witten invariants of $\overline{{\rm Graph} (\mu)}_{t}$
that we can explicitly compute.

\subsection{The neighborhood}\text{} 
The transversal $A_3$-singularity is identified with $\IP^1$
 by the morphism 
$[z_0:z_1]\mapsto [0:0:z_0:z_1]$.
By abuse of notation we denote by the same symbol $\cO(a)$
the sheaf $\cO_{\IP^1}(a)$ and the corresponding vector bundle, for any $a\in \IZ$.
Moreover we identify $\mathcal{O}(a)\otimes\mathcal{O}(b)$ with 
$\mathcal{O}(a+b)$ using the canonical isomorphism.
For any vector bundle $E$, we denote by $\ul{0}_{E}$ its zero section.
Finally, we set 
\begin{displaymath} 
  V_i:=\{[x_0:x_1:x_2:x_3]\in |\IP(1,3,4,4)| \mbox{ such that } x_i\not=0 \}  
\end{displaymath} 
for any $i\in \{0,1,2,3\}$, and 
$$
V:=V_2\cup V_3 \subset |\IP(1,3,4,4)|. 
$$

Consider the bundle morphism 
\begin{align*} 
  \psi:\mathcal{O}(1)\oplus\mathcal{O}(3)\oplus \mathcal{O}(1) 
  &\longrightarrow \mathcal{O}(4) \\(\xi,\eta,\zeta) &\longmapsto \xi 
  \ot \eta - \zeta^{\ot 4},
\end{align*} 
and the inverse image under $\psi$ of the zero section of $\cO(4)$:
$\psi^{-1}(\rm{Im}(\ul{0}_{\cO(4)} ))\subset \mathcal{O}(1)\oplus\mathcal{O}(3)\oplus \mathcal{O}(1)$. 
We have the following
\begin{lemma}\label{lem:iso} 
  The variety $V$ is isomorphic to $\psi^{-1}(\rm{Im}(\ul{0}_{\cO(4)} ))$. 
\end{lemma} 
\begin{proof} 
  An easy computation shows that $V_{2}\simeq\Spec\left( 
    \IC[s,u,v,w]/{(uv-w^4)} \right)$ and $V_3 \simeq 
  \Spec \left( \IC[t,x,y,z]/(xy-z^4)\right)$. 
The affine open subvarieties $V_2,V_3\subset V$ 
glue together by means of the following ring isomorphism 
\begin{align*}  
 \fr{\IC[s,\fr{1}{s},u,v,w]}{(uv-w^4)}  &\longrightarrow    
 \fr{\IC[t,\fr{1}{t},x,y,z]}{(xy-z^4)} \\  
s& \longmapsto  \fr{1}{t}\\ 
u& \longmapsto  \fr{1}{t}x\\ 
v& \longmapsto  \fr{1}{t^3}y\\ 
w& \longmapsto  \fr{1}{t}z.  
\end{align*}  
On the other hand, consider a trivialization of the bundle 
$\mathcal{O}(1)\oplus\mathcal{O}(3)\oplus \mathcal{O}(1)$ on 
$W_{0}=\{[z_{0},z_{1}]\in \IP^1\mid z_{0}\neq 0\}$.
On such a 
trivialization, the morphism $\psi$ is given by
\begin{align*} 
  W_{0}\times\IC^{3} &\longrightarrow W_{0}\times\IC\\ 
(s,v_{1},v_{2},v_{3})&\longmapsto (s,v_{1}v_{2}-v_{3}^{4}). 
\end{align*} 
Hence we have that, over $W_0$, 
$\psi^{-1}(\rm{Im}(\ul{0}_{\mathcal{O}(4)}))$ is 
$\Spec \left( \IC[s,v_{1},v_{2},v_{3}]/(v_{1}v_{2}-v_{3}^{4}) \right)$.  
If we do the same over 
$W_{1}=\{[z_{0},z_{1}]\in \IP^1\mid z_{1}\neq 0\}$,  we deduce 
that $V$ and $\psi^{-1}(\ul{0}_{\mathcal{O}(4)})$ are 
union of the same affine varieties with the same gluing. This proves 
 that they are isomorphic. 
\end{proof} 

\subsection{The deformation} \text{} 
 
We construct a deformation of $V$.
The construction is inspired by the theory of deformations
of surfaces with  ${\rm A}_n$ singularities, for this
our reference is \cite{Ty}. 

We define 
\begin{equation}\label{fibrazione}
f:V\ra \IP^1
\end{equation}
to be the composition
of the isomorphism $V\xrightarrow{\cong}\psi^{-1}(\rm{Im}(\ul{0}_{\cO(4)} ))$
in Lemma \ref{lem:iso}, followed by the inclusion
$\psi^{-1}(\rm{Im}(\ul{0}_{\cO(4)} ))\subset \cO(1)\op \cO(3)\op \cO(1)$,
and then the bundle map. The following remark is crucial.
\begin{remark}\nf
The morphism \eqref{fibrazione} exhibits $V$
as a $3$-fold fibered over $\IP^1$ with fibers isomorphic to
a surface ${\rm A}_3$ singularity, furthermore the 
fibration is locally trivial.
\end{remark} 

The aim is to extend some of the results of \cite{Ty}  to $V$, 
when viewed
as a family of such surfaces with respect to \eqref{fibrazione}.

Consider the bundle morphism 
\begin{align*}
\chi:\mathcal{O}(1)^{\oplus 4} &\longrightarrow 
  \mathcal{O}(1) \\ 
  (\delta_{1},\ldots 
  ,\delta_{4})&\longmapsto\delta_{1}+\cdots+\delta_{4},
  \end{align*} 
and set  $\mathcal{F}:=\chi^{-1}(\rm{Im}(\ul{0}_{\mathcal{O}(1)}))$. 
Then consider the bundle morphism 
\begin{align*}
\pi:\mathcal{O}(1)\oplus\mathcal{O}(3)\oplus\mathcal{O}(1)\oplus 
  \mathcal{F}&\longrightarrow 
  \mathcal{O}(4) \\ 
  (\xi,\eta,\zeta,\delta_{1},\ldots 
  ,\delta_{4})&\longmapsto\xi\otimes\eta-\otimes_{i=1}^{4}(\zeta+\delta_{i})
  \end{align*} 
and set
  $\mathcal{V}_{\mathcal{F}}:=\pi^{-1}(\rm{Im}(\ul{0}_{\mathcal{O}(4)}))$. 
We obtain the following Cartesian diagram 
\begin{equation*} 
\begin{CD}
  V @>>> {\cV}_{\cF} \\
  @VVfV @VVFV \\
  \IP^1 @>\ul{0}_{\cF}>> \cF
\end{CD}
\end{equation*} 
where $f$ is defined in \eqref{fibrazione}
and the vertical right hand side arrow is the composition of the inclusion
 ${\cV}_{\cF}\ra \cO(1)\op \cO(3)\op \cO(1)\op \cF$
followed by the projection $\cO(1)\op \cO(3)\op \cO(1)\op \cF \ra \cF$.
 
\medskip
Note that $F:\cV_{\cF}\ra \cF$ is a family of surfaces.
We now construct a simultaneous resolution. 
Consider the rational map 
\begin{align*} 
  \mu : \mcal{V}_{\mcal{F}} &\dashrightarrow \IP (\mcal{O}(1)\op 
  \mcal{O}(1))\times \IP (\mcal{O}(1)\op \mcal{O}(2))\times 
  \IP (\mcal{O}(1)\op \mcal{O}(3))\\ 
  (\xi,\eta,\zeta,\de_1,...,\de_{4}) &\longmapsto (\xi,\zeta 
  +\de_1)\times (\xi,(\zeta +\de_1)\ot (\zeta +\de_2)) \times (\xi, 
  \ot_{i=1}^3(\zeta +\de_i)), 
\end{align*} 
and let $\Graph(\mu)$ be the graph of $\mu$. \\Then take the closure of 
$\Graph(\mu)$ in $\mcal{V}_{\mcal{F}}\times \left( \times_{i=1}^3 
\IP (\mcal{O}(1)\op \mcal{O}(i)) \right)$: 
$$ 
\overline{\Graph(\mu)} \subset \mcal{V}_{\mcal{F}}\times \left( \times_{i=1}^3 
\IP (\mcal{O}(1)\op \mcal{O}(i))\right). 
$$ 
Let $\cR : \overline{{\rm Graph}(\mu)} \ra \cV_\cF$ be the composition
of the inclusion \\$\overline{{\rm Graph}(\mu)} \ra \cV_\cF \times 
\left( \times_{i=1}^3 \IP(\cO(1) \op \cO(i))\right)$
followed by the projection on the first factor
$\cV_\cF \times \left( \times_{i=1}^3 \IP(\cO(1) \op \cO(i))\right) \ra \cV_\cF$.
\begin{lemma}
The following diagram
\begin{equation} 
  \label{sr,bis} 
  \xymatrix{\overline{\Graph(\mu)} \ar[r]^>>>>>{\cR}\ar[d]^{F\circ \cR}& \ar[d]^{F}\mathcal{V}_{\mathcal{F}} \\  
\mathcal{F} \ar[r]^{\mathrm{id}}& \mathcal{F}} 
\end{equation} 
is a simultaneous resolution  of $F:\mcal{V}_{\mcal{F}}\ra \mcal{F}$.
\end{lemma}
\noindent \textbf{Proof.} The property of being a simultaneous 
resolution is local in $\mcal{F}$. The diagram \eqref{sr,bis}
is fibered over $\IP^1$. If we restrict it to an open subset of $\IP^1$ where 
$\mcal{O}(1)$ is trivial, then the assertion is exactly the result of 
E. Brieskorn \cite{B}. \qed 

\bigskip
Set $\Delta:=\IC$. For any section $\theta \in H^0(\IP^1, \mcal{F})$, we get a deformation 
of $V$ parametrized by $\De$ as follows 
\begin{displaymath} 
  \xymatrix{V \ar[r]\ar[d]^{f}&\mathcal{V}_{\theta}\ar[r]\ar[d]^{f_{\theta}} & \mathcal{V}_{\mathcal{F}}\ar[d]^F \\ 
\IP^{1} \ar[r]& \IP^{1}\times \Delta\ar[r]^-{\Theta} & \mathcal{F}} 
\end{displaymath} 
where $\Theta :\IP^1 \times \De \ra \mcal{F}$ sends $([z_0:z_1],t)$ to $t\cd 
\theta([z_0:z_1])$, and $\mcal{V}_{\theta}$ is defined by the requirement that the diagram is 
Cartesian, the map $\IP^1 \ra \IP^1\times \De$ is the inclusion
$[z_0:z_1]\mapsto ([z_0:z_1],0)$.
The pull-back of the diagram \eqref{sr,bis} with respect to $\Theta$ gives the 
following diagram  
\begin{equation}\label{d,bis}  
\begin{CD}  
\overline{\Graph(\mu)}_{\theta} @>\rho_{\theta}>> \mcal{V}_{\theta}\\  
@VVV @VV f_\theta V \\  
\IP^1 \times \De @>\rm{id}>> \IP^1 \times \De  
\end{CD}  
\end{equation}  
where $\rho_{\theta}$ is the pull-back of $\cR$ in \eqref{sr,bis}.
We remark that \eqref{d,bis} is a simultaneous resolution of $\mcal{V}_{\theta}$ over $\IP^{1}\times\Delta$.
  
\subsection{Computation of the invariants} \text{} 

We specialize the previous construction in the case where 
$\theta$ is given as follows.
Let $\de \in H^0(\IP^1, \mcal{O}(1))$ be a nonzero section, and set
\begin{displaymath} 
  \de_\ell:=\rm{exp}\left( \fr{(2\ell+1)\pi i}{4} \right)\cd \de, 
\quad \ell\in \{ 1,...,4 \}, 
\end{displaymath} 
then set $\theta:=(\de_1,...,\de_{4})\in H^0(\IP^1,\mcal{F})$. 
For any $t\in \De$, let us denote $\cV_t:= f^{-1}_\theta\left(\IP^1\times \{t\}\right)$,
$f_t:\cV_t \ra  \IP^1\times\{t\}$ the restriction of $f_\theta$,
$\overline{{\rm Graph}(\mu)}_t:= \rho_\theta^{-1}(\cV_t)$
and $\rho_t:\overline{{\rm Graph}(\mu)}_t\ra \cV_t$ the restriction of $\rho_\theta$.
We have the following commutative diagram 
\begin{displaymath}\label{above}
  \xymatrix{\rho^{-1}(V)=\overline{\Graph(\mu)}_{0}  \ar@{^{(}->}[r] \ar[d]^{\rho_{0}=\rho}& 
    \overline{\Graph(\mu)}_{\theta} \ar[d]^{\rho_{\theta}}&  
\ar@{_{(}->}[l] \ar[d]^{\rho_{t}} \overline{\Graph(\mu)}_{t} \\ 
V \ar@{^{(}->}[r]\ar[d]^{f_0=f}& \ar[d]^{f_\theta}\mathcal{V}_{\si}&\ar[d]^{f_t}\ar@{_{(}->}[l] \mathcal{V}_{t} \\ 
\IP^{1}\times\{0\} \ar@{^{(}->}[r]& \IP^{1}\times \Delta &\ar@{_{(}->}[l] \IP^{1}\times \{t\}   
} 
\end{displaymath} 
 
\begin{lemma}\label{lem:graph,t} 
  Let $\delta$ be a global section of $\mathcal{O}(1)\to \IP^{1}$ that 
  vanishes only at one point. Then, for $t\neq 0$, the variety 
  $\overline{\Graph(\mu)}_t$ has only one connected nodal complete 
  curve of genus $0$ whose dual graph is of type $A_{3}$ and which is contracted by 
  $\rho_{t}$ (see the diagram above). 
\end{lemma} 
\begin{proof}
  Without lost of generality we can assume that $\delta$ vanishes 
  only at the point $[1:0]$.  Let $W_{0}:=\{[z_{0}:z_{1}]\in 
  \IP^{1}\mid z_{0}\neq 0\}$. As our bundles are trivial over $W_{0}$, 
  the restriction of  $V$ over $W_0$ is given by $W_{0}\times \IV(xy-z^{4}) \subset 
  W_{0}\times\IC^{3}$.  The choice of the $\delta_{\ell}$ implies 
  that the $3$-fold $\mathcal{V}_{t}$ is given by 
\begin{displaymath} 
  W_{0}\times\IV(xy-\prod_{\ell=1}^{4}(z+\delta_{\ell}t))= 
  W_{0}\times\IV(xy-z^{4}-(t\delta)^{4}) \subset W_{0}\times\IC^{3}. 
\end{displaymath} 

By means of $f_t$, $\mcal{V}_t$ is viewed a family of 
surfaces parametrized by $\IP^1$. As $t\neq 0$ and $\delta([1:0])=0$, the only 
singular surface of the family is the surface 
$f^{-1}_{t}([1:0]\times \{t\})$, which is a surface with an isolated
$A_{3}$-singularity. As $\rho_{t}:\overline{\Graph(\mu)}_{t}\to \mathcal{V}_{t}$ is a 
simultaneous resolution over $\IP^{1}\times\{t\}$, the fiber 
$\overline{\Graph(\mu)}_{([1:0],t)}$ is a smooth surface with only one 
complete connected curve of genus $0$ whose dual graph is of type $A_{3}$ and 
which is contracted by $\rho_{t}$. For 
any $[z_{0}:z_{1}]\neq [1:0]$, the fiber 
$\overline{\Graph(\mu)}_{([z_{0}:z_{1}],t)}$ is isomorphic to the 
smooth surface $f^{-1}_{t}([z_{0}:z_{1}]\times\{t\})$.  Hence, the 
exceptional locus of the resolution $\rho_{t}:\overline{\Graph(\mu)}_t 
\to\mathcal{V}_{t}$ has  only one connected nodal complete curve of 
genus $0$ whose dual graph is of type $A_{3}$ and which is contracted 
by $\rho_{t}$. 
\end{proof} 

Let $\Ga_1,\Ga_2,\Ga_3 \in H_2(\overline{{\rm Graph}(\mu)}_t; \IZ)$
be the homology classes of the components 
of the connected nodal complete curve of genus
zero whose dual graph is of type ${\rm A}_3$ and which is contracted by
$\rho_t$. Let us assume that they are numbered in such a way that,
if $\Ga_i$ is the class of $\ti{\Ga}_i$, then the  intersection $\ti{\Ga}_i\cap \ti{\Ga}_j$
is empty if $\mid i-j\mid >1$.
Then the previous Lemma implies that, for $t\neq 0$, $\overline{\Graph(\mu)}_t$ satisfies the 
hypothesis of Proposition 2.10 of \cite{BKL}. Therefore we deduce the following formula:  
\begin{align}\label{eq:deg,virt} \rm{deg} [\overline{\mcal{M}}_{0,0}(\overline{\Graph(\mu)}_t,\Ga)]^{\rm{vir}}=  
\begin{cases}  
  {1}/{d^3} &  \mb{if} \, \Ga=d(\Gamma_{\mu} +\Ga_{\mu+1}+...+\Ga_{\nu}), \, \mb{for} \, \mu\leq \nu;\\ 
  0 &  \mb{otherwise}. 
\end{cases}  
\end{align} 
This formula together with the next lemma completes the proof of Theorem \ref{GW1344}. 
 
\begin{lemma} 
For any $t\in \De$, the following equality holds 
$$ 
\rm{deg}[\overline{\mcal{M}}_{0,0}(\overline{\rm{Graph}(\mu)}_t,\Ga)]^{\rm{vir}}= 
\rm{deg}[\overline{\mathcal{M}}_{0,0}(Z,\Ga)]^{\rm{vir}}. 
$$ 
\end{lemma}   
\begin{proof} 
  Since $\Ga$ is the homology class of a contracted curve, 
we have an isomorphism of moduli stacks (see  
\cite[Lemma $7.1$]{P}):
\begin{equation} \label{ims}
  \overline{\mcal{M}}_{0,0}(Z,\Ga) \simeq \overline{\mcal{M}}_{0,0}(\rho^{-1}(V),\Ga);
\end{equation} 
in particular the right hand side moduli stack is proper with  
projective coarse moduli space. The isomorphism \eqref{ims}
identifies the tangent-obstruction theories used to define the 
Gromov-Witten invariants,
hence the virtual fundamental classes $[\overline{\mcal{M}}_{0,0}(\rho^{-1}(V),\Ga)]^{\rm{vir}}$ 
and $[\overline{\mcal{M}}_{0,0}(Z,\Ga)]^{\rm{vir}}$ have the same degree. 
Then it is enough to prove that, for any $t\in \Delta$,
\begin{equation}\label{eq:4} 
\rm{deg}[\overline{\mcal{M}}_{0,0}(\overline{\rm{Graph}(\mu)}_t,\Ga)]^{\rm{vir}}= 
\rm{deg}[\overline{\mcal{M}}_{0,0}(\rho^{-1}(V),\Ga)]^{\rm{vir}}. 
\end{equation} 
 
Gromov-Witten invariants of 
projective varieties are invariant under deformation of the target 
variety. We now explain why this result holds for $\rho^{-1}(V)$ and $\overline{\Graph(\mu)}_{t}$ 
even if they are not projective.
 
Let $q_{\theta}:\overline{\rm{Graph}(\mu)}_{\theta} \ra \De$
be the composition of $f_\theta \circ \rho_\theta$ in \eqref{d,bis} 
followed by the projection $\IP^1 \times \De \ra \De$.
The morphism $q_\theta$ is smooth as composition  of  smooth morphisms.
Moreover $q_{\theta}$ factories through an embedding followed by a 
projective morphism. To see this, it is enough to prove the same statement
for the morphism $F\circ \cR: \overline{\Graph(\mu)}\to \mathcal{F}$ in \eqref{sr,bis}.
By construction, $\overline{{\rm Graph}(\mu)}$
is embedded in $\cO(1)\op \cO(3)\op \cO(1)\op \cF \times \times_{i=1}^3 \IP(\cO(1)\op \cO(i))$,
moreover $F\circ \cR$ is the restriction of the projection 
$\cO(1)\op \cO(3)\op \cO(1)\op \cF \times 
\left( \times_{i=1}^3 \IP(\cO(1)\op \cO(i)) \right) \ra \cF$.\\
Let us consider now the projection $\cO(1)\op \cO(3)\op \cO(1)\op \cF \ra \cF$,
it has  a vector bundle structure over $\cF$, then it can be seen
as a subbundle of the projective bundle
$\IP(\cO(1)\op \cO(3)\op \cO(1)\op \cF \op \cO_\cF) \ra \cF$,
therefore we have that $F\circ \cR$
factors as the composition of an embedding followed by 
a projective morphism.
 
To finish the proof, let us consider the moduli stack which 
parameterizes relative stable maps to $q_\theta:\overline{\Graph(\mu)}_{\theta}\ra \De$ of 
homology class $\Ga$ and genus zero.  We denote it by 
$\overline{\mcal{M}}_{0,0}(\overline{\Graph(\mu)}_{\theta}/\De,\Ga)$.  As 
$\Ga$ is the class of curves which are contracted by the resolution 
$\rho_{\theta}$ and $q_{\theta}:\overline{\Graph(\mu)}_{\theta}\to 
\Delta$ factories through an embedding followed by a projective 
morphism, Theorem $1.4.1$ of \cite{AVcssm} implies that the moduli 
space 
$\overline{\mcal{M}}_{0,0}(\overline{\Graph(\mu)}_{\theta}/\De,\Ga)$ is a 
proper Deligne-Mumford stack. Since the class $\Ga$ is contracted by 
$\rho_{\theta}$, for any $t \in \Delta$ the fiber at $t$ of the natural 
morphism 
$\overline{\mcal{M}}_{0,0}(\overline{\Graph(\mu)}_{\theta}/\De,\Ga)\ra 
\De$ is the proper Deligne-Mumford stack 
$\overline{\mcal{M}}_{0,0}(\overline{\Graph(\mu)}_{t},\Ga)$. 
Then the same proof of Theorem 4.2 in \cite{LTvmc}
applies in this situation and we get \eqref{eq:4}. 
\end{proof} 
   
 \section{ The case $\cX= \IP(1,\ldots,1,n)$}\label{sec:p11n} 
In this Section we will prove the following proposition.  
Here $\cX$ denotes $\IP(1,...,1,n)$.
\begin{proposition}\label{p11n}
Let $n\geq 2$ be an integer and consider the $n$-dimensional weighted projective space
$\IP(1,...,1,n)$. Let $Z$ be the 
crepant resolution of $|\IP(1,...,1,n)|$ defined in point (1) below.    
Then, there is a ring isomorphism  
$$ 
{\rm H}^\star(Z;\IC)\cong {\rm H}^\star_{\rm CR}(\IP(1,\ldots,1,n);\IC). 
$$ 
\end{proposition} 
 
\begin{proof}We follow the steps described in Section \ref{sec:CCRC}. 
   
 The coarse moduli space $|\IP(1,\ldots,1,n)|$ has an  
isolated singularity of type $\frac{1}{n}(1,\ldots,1)$ at the point $[0:\ldots:0:1]$.  
 
\medskip 
 (1) We identify the stacky fan $(N,\be,\Si)$ defined in \eqref{stackyfan} with 
$(\IZ^n,\{\Lam(\be(v_i))\}_{i=0}^n , \Si)$ by means of the isomorphism $\Lam :N \ra \IZ^n$
defined  by sending $v_0$ to $(-1,...,-1,-n)$ and $v_i$ to the 
$i$-th vector of the standard basis of $\IZ^n$,
for $i\in \{ 1,...,n \}$.

The crepant resolution is defined as follows:
consider the ray $\lan P\ran$ generated by 
${\rm P}:=(0,...,0,-1)=\fr{1}{n}\sum_{i=0}^{n-1}\Lam(\be(v_i))$,
then let $\Si'$ be the fan obtained from $\Sigma$ by replacing the cone 
generated by $\Lam(\be(v_0)),\ldots,\Lam(\be(v_{n-1}))$ with the cones 
generated by $\Lam(\be(v_0)),\ldots,\widehat{\Lam(\be(v_i))},\ldots,\Lam(\be(v_{n-1}))$ and $\rm P$ for  
any $i\in \{ 0,\ldots,n-1 \}$. We draw as an example the polytope 
for the case $n=3$ in Figure \ref{fig5}. 
Define $Z$ to be the toric variety associated to  $\Si'$,
and $\rho:Z\ra X$ to be the morphism associated to the identity in $\IZ^n$.  
  
 \begin{figure}[!ht]  
 \begin{center}  
 \psfrag{A}{$\Lam(\be(v_{1}))$} \psfrag{B}{$\Lam(\be(v_{2}))$} 
\psfrag{C}{$\Lam(\be(v_{3}))$} \psfrag{D}{$\Lam(\be(v_{0}))$}  
 \psfrag{E}{$P$}  
  \includegraphics[width=0.4\linewidth]{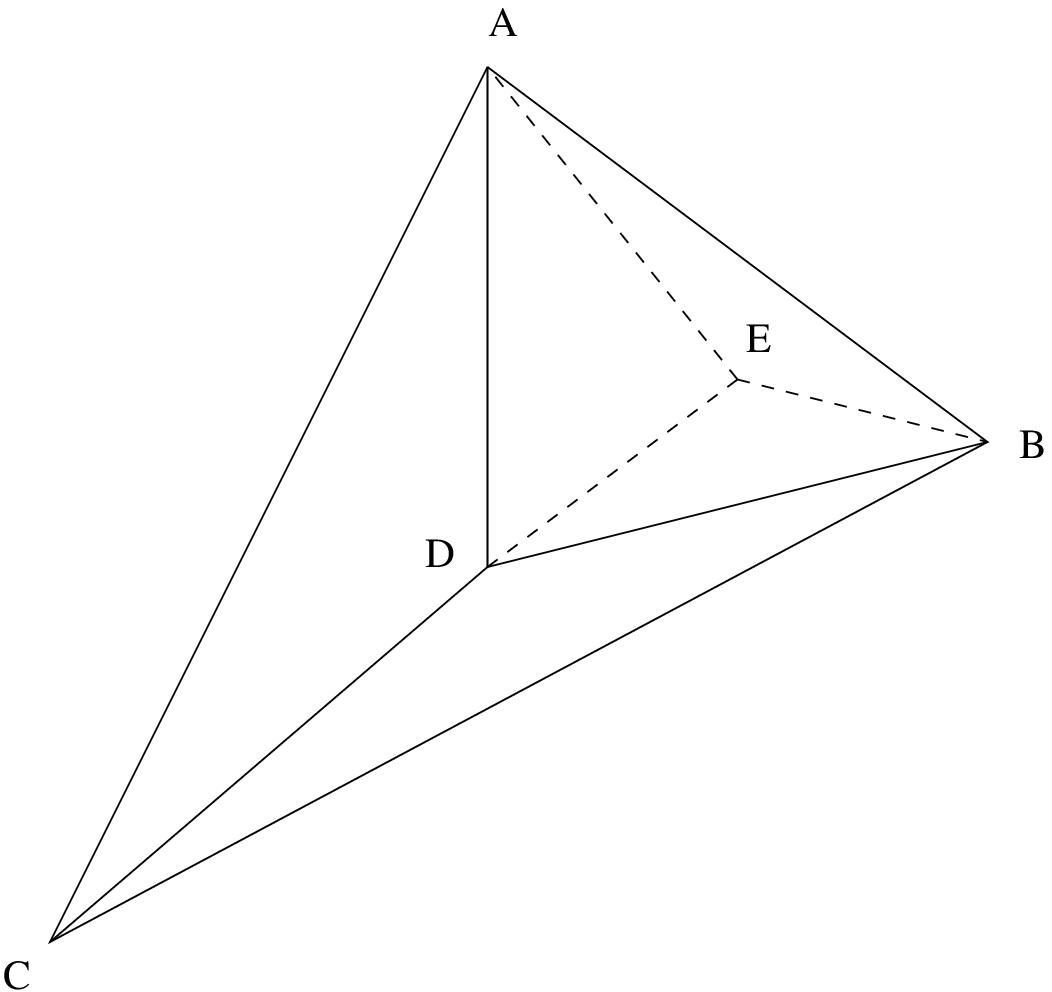}  
 \end{center}\caption{Polytope of $\IP(1,1,1,3)$ and a crepant resolution}\label{fig5}  
 \end{figure}  
  
\medskip
 (2) Let $b_i\in H^2(Z;\IC)$ ($e\in H^2(Z;\IC)$ resp.)
be the first Chern class of the line bundle associated to 
the torus invariant divisor corresponding to the ray generated by 
$\Lam(\be(v_i))$ ($\rm P$ resp.) for any $i\in\{0,...,n\}$.
We have $H^\star(Z;\IC)\cong \IC[b_0,\ldots,b_n,e]/I$ where $I$ is generated by:  
 \begin{eqnarray*}  
 -b_0+b_i\quad \text{ for } 1\leq i\leq n-1,\\  
 -nb_0-e+b_n, \, eb_n, \, b_0\cdots b_{n-1}.  
 \end{eqnarray*}  
Set $h:=\frac{1}{2n}(b_0+\ldots+b_n+e)=b_0+\frac{1}{n}e$, then we get:  
 $$  
 H^\star(Z;\IC)\cong \IC[h,e]/\langle  
 h^n+(-1)^n\left(\frac{e}{n}\right)^n,he\rangle.  
 $$  

\medskip  
 (3) ${\rm M}_{\rho}(Z)$ is generated by one class
$\Ga_1:=\PD\left(\left(h-\frac{e}{n}\right)^{n-2}e\right)$.  

\medskip  
 (4) We will set the quantum parameter $q_1=0$, then we do not have
to compute any non trivial Gromov-Witten invariant.

\medskip
 (5) 
We follow the description given in \cite{BMPmodel} of the 
Chen-Ruan cohomology ring. The twisted sectors are indexed by the set
${\rm T}=\left\{ {\rm exp}\left(\fr{2\pi \ii k}{n}\right) 
\, \mid \, k\in\{0,...,n-1\}\right\}$.
For any $g\in {\rm T}-\{ 1 \}$, $\cX_{(g)}\cong \IP(n)$,
while $\cX_{(1)}\cong \cX$.
As vector space we have
\begin{equation}\label{crcohomologyp11n}
H^\star_{\rm CR}(\cX;\IC):= \op_{g\in T}H^\star(\cX_{(g)}).
\end{equation}
Let 
$$
H,\, E_1 \in H^2_{\rm CR}(\cX;\IC)
$$
be the image of $c_1\left(\cO_{\cX}(1)\right)\in H^2(\cX;\IC)$,
$1\in H^0(\cX_{\left({\rm exp}\left(\fr{2\pi \ii }{n} \right) \right)};\IC)$
respectively with respect to the inclusion
$H^\star(\cX_{(g)})\ra H^\star_{\rm CR}(\cX)$
determined by \eqref{crcohomologyp11n}.
Then we have the following presentation:
 $$  
 H^\star_{\CR}(\IP(1,\ldots,1,n);\IC)\cong 
\IC\left[H,E_1\right]/\left\langle H^n-
\left(E_1\right)^n,HE_1\right\rangle.  
 $$  
 
\medskip 
 (6) The ring isomorphism  
 $$  
 H^\star_{\CR}(\IP(1,\ldots,1,n);\IC)\xrightarrow{\sim}  
 H^\star(Z;\IC)  
 $$  
 is obtained by mapping $H\mapsto h$ and
$E_1\mapsto  -\exp\left(\frac{\ii\pi}{n}\right)\frac{e}{n}$.  
\end{proof}

\bibliographystyle{amsalpha}  
\bibliography{BMPBib} 
  
\end{document}